\documentclass{gtpart}

\newcommand{\heuteIst}{Oct 16, 2000 }

\newcommand{\longversion}[1]{{#1}}

\usepackage{amsmath,amsthm,amsfonts,latexsym,amscd,amssymb,enumerate}
\usepackage{amsthm}

\title{Finite group extensions and the Baum-Connes conjecture}
\author{
Thomas Schick}
\givenname{Thomas}
\surname{Schick}
\address{Mathematisches Institut\\ Georg-August-Universit\"at G\"ottingen\\
  Bunsenstr.~3--5\\ 37073 G\"ottingen\\ Germany}

\email{schick@uni-math.gwdg.de}
\urladdr{http://www.uni-math.gwdg.de/schick/}

\arxivreference{math.KT/0209165}
\arxivpassword{westx}
        
\swapnumbers
\theoremstyle{plain}
\newtheorem{theorem}{Theorem}

\newtheorem{corollary}[theorem]{Corollary}
\newtheorem{proposition}[theorem]{Proposition}

\theoremstyle{definition}
\newtheorem{definition}[theorem]{Definition}
\newtheorem{example}[theorem]{Example}

\theoremstyle{remark}
\newtheorem{remark}[theorem]{Remark}

\usepackage{amsmath,amsfonts,latexsym}

\newcommand{\BCclass}{\mathbf{LH}\mathcal{ETH}}

\newcommand{\complexs}{\mathbb{C}}

\newcommand{\integers}{\mathbb{Z}}

\newcommand{\generate}[1]{\langle#1\rangle}

\newcommand{\subgroup}{<}

\newcommand{\semiProd}{\rtimes}
\newcommand{\semiprod}{\semiProd}

\newcommand{\forget}[1]{}

\numberwithin{equation}{section}

\newcommand{\extendableGroups}{\mathcal{F}}
\newcommand{\GoodextendableGroups}{\mathcal{F}^+}

\begin{document}
\date{Last compiled \today; last edited  \heuteIst or later}

\begin{abstract}
    In this note, we exhibit a method to prove the Baum-Connes
    conjecture (with coefficients) for extensions with finite
    quotients of certain groups which already satisfy the Baum-Connes
    conjecture. Interesting examples to which this method
    applies are torsion-free finite extensions of the pure braid
    groups, e.g.~the full braid groups, and fundamental groups of certain link
    complements in $S^3$.
\end{abstract}

\maketitle

The Baum-Connes conjecture 
for a group
$G$ states that the Baum-Connes map
\begin{equation}\label{eq:Baum_Connes}
  \mu_r\colon K^G_*(\underline{E}G,A) \to K_*(C^*_{red}(G,A))
\end{equation}
is an isomorphism for every $C^*$-algebra $A$ with an action of $G$ by
$C^*$-algebra homomorphisms. In this note the term ``Baum-Connes conjecture'' will always mean
``Baum-Connes conjecture with coefficients'', and all group are assumed to be
discrete and countable.  

\longversion{Here, the left hand side is the
equivariant K-homology with coefficients in $A$ of the universal
space $\underline{E} G$
for proper $G$-actions, which is homological in nature. The right hand
side, the K-theory of the reduced crossed product of $A$ and $G$,
belongs to the world of $C^*$-algebras
and ---to some extend--- representations of groups. If $A=\complexs$
with the trivial action, the right hand side becomes the K-theory of
the reduced $C^*$-algebra of $G$. If, in addition, $G$ is
torsion-free, the left hand side is the K-homology of the
classifying space of $G$.

The Baum-Connes conjecture has many important connections to other
questions and areas of mathematics. The injectivity of $\mu_r$ implies the
Novikov conjecture about homotopy invariance of higher signatures. It
also implies the stable Gromov-Lawson-Rosenberg conjecture about the
existence of metrics with positive scalar curvature on
spin-manifolds. The surjectivity, on the other hand, gives information
in particular about $C^*_{red}G$. If $G$ is torsion-free, it implies
e.g.~that this $C^*$-algebra contains no
idempotents different from zero and one. Since we are only considering
the Baum-Connes conjecture with coefficients, all these properties
follow for all subgroups of $G$, as well.
}

We do not want to repeat the construction of the $K$-groups and the map in
the  Baum-Connes conjecture \eqref{eq:Baum_Connes}, instead, the reader
is referred to \cite{Baum-Connes-Higson(1993),Julg(1998)}.  
\longversion{
Higson, Lafforgue and Skandalis \cite{MR1911663}, using groups
constructed by Gromov \cite{MR1978492}, have produced
counterexamples to the conjecture (with non-trivial, commutative
coefficients). However, we will
concentrate on groups for which the Baum-Connes conjecture is known to
be true and will prove it for new examples.
}
\longversion{A well known facts is
that the Baum-Connes conjecture  is inherited by
arbitrary subgroups, and there are rather precise results describing
its behavior under
group extensions (we will recall this below). However, it is
almost completely unknown what is
happening for an extension with finite quotient of a group which
satisfies the Baum-Connes conjecture. The most
prominent open examples are probably the full braid groups, which
contain the pure braid groups as subgroups of finite index. For the
latter, the Baum-Connes conjecture is well known to be true.
}

The main goal of this note is to prove Baum-Connes 
for the full braid groups, and for other classes of groups which arise
as (finite) extensions of groups for which Baum-Connes is known. We
use the following results:

\begin{theorem}\label{theo:subgroups}
  Assume $G$ is a  group which satisfies the
  Baum-Connes conjecture with coefficients. If $H$ is a subgroup of
  $G$, then $H$ satisfies the Baum-Connes conjecture with coefficients, too.
\end{theorem}

 \begin{theorem}\label{theo:trees}\cite{Oyono-Oyono(1998)}
   Assume $G$ is a group acting on a tree. The Baum-Connes
   conjecture  is true for $G$ if and only if it is
   true for every isotropy subgroup of the action on the vertices.
 \end{theorem}
Note that, from a logical point of view, Theorem \ref{theo:subgroups} is a
consequence of Theorem \ref{theo:trees}.

\begin{theorem}\label{theo:extensions}(\cite[Section
  3]{Chabert-Echterhoff(2000)},\cite[Theorem
  3.1]{Oyono-Oyono(2000)} 
  Assume we have an extension of groups 
  \begin{equation*}
1\to H\to G\xrightarrow{\pi} Q\to
  1.
\end{equation*}
  For every finite subgroup $E$ of $Q$ let $H_E\subgroup G$ be the
  inverse image under $\pi$ of $E$ in $G$. Assume $Q$ and all groups
  $H_E$ satisfy the
  Baum-Connes conjecture with coefficients.

  Then the Baum-Connes conjecture with coefficients is also true for $G$.
\end{theorem}

\begin{theorem}\label{theo:amenable}\cite{Higson-Kasparov(1997),Higson-Kasparov(1997a)}
  (Compare also \cite{Julg(1998)}):
  Let $G$ be an amenable group. Then the Baum-Connes
  conjecture with coefficients holds for $G$.
\end{theorem}

Theorem \ref{theo:extensions} and Theorem \ref{theo:amenable} together
immediately imply the following corollary.
\begin{corollary}\label{corol:tfamext}
  Assume we have an extension $1\to H\to G\to A\to 1$ of
  groups, $A$ is torsion-free and amenable, and the
  Baum-Connes conjecture with coefficients  holds for $H$. Then it also
  holds for $G$.
\end{corollary}

Our (naive) idea to deal with finite group extensions is to
exhibit them as extension with torsion-free quotient. More precisely,
assume we have an extension $1\to H\to G\to Q\to 1$ where $H$
satisfies the Baum-Connes conjecture  and $Q$ is a
finite group. Our goal is, to find a normal subgroup $U$ of $G$ which
is contained in $H$ and such that $G/U$ is torsion-free and
amenable\longversion{ (e.g.~virtually nilpotent or virtually solvable)}. 
If $H$ is torsion-free, a necessary condition for the existence of
such a factorization is
that $G$ is torsion-free, too. We will observe that there are large classes of
torsion-free groups such that every torsion-free finite extension
admits a factorization of the type we are looking for, including all
pure braid groups.

 We need the following
notation:
\begin{definition}\label{def:complete}
  Let $H$ be a group and $p$ a prime number. Let ${\hat H}^p$
  denote the pro-$p$ completion of $H$, i.e.~the inverse limit of the
  system of
  finite $p$-group quotients of $H$. There is a natural homomorphism $H\to
  {\hat H}^p$\longversion{ (not necessarily injective)}.
  The Galois cohomology
  $H^*({\hat H}^p,\integers/p)$ is defined to be the direct limit of
  the cohomology of the finite $p$-group quotients of $H$ (with
  $\integers/p$-coefficients). For the theory of profinite groups and
  Galois cohomology compare e.g.~\cite{Wilson(1998)}.

  The group $H$ is called \emph{cohomologically complete} if the
  natural homomorphisms $H^*({\hat H}^p,\integers/p)\to
  H^*(H,\integers/p)$ are isomorphisms for every prime $p$.
\end{definition}

The second property of $H$ we are going to use is the existence of
many quotients which are torsion-free and amenable.
\begin{definition}\label{def:enough_quotients}
  Let $H$ be a group. We say $H$ has \emph{enough
  amenable torsion-free quotients} if for every normal subgroup $U$ of $H$ of
  finite $p$-power index  (for some prime $p$) another normal subgroup
  $V$ of $H$ exists which is contained in $U$ and such that $H/V$ is
  torsion-free and elementary amenable. We say that $H$ has \emph{enough
    nilpotent torsion-free quotients} if we find $V$
  such that $H/V$ is nilpotent.
\end{definition}

\longversion{
\begin{remark}
  Recall that the class of elementary amenable groups is the smallest
  class of groups which contains all finite and all abelian groups and
  which is closed under extensions and directed unions. It contains in
  particular all nilpotent and all solvable-by-finite groups.

  For the purpose of our paper, we could weaken the condition to $H/V$
  being amenable. However, this would deviate from the notation used
  in \cite{Linnell-Schick(2000)}. Moreover, all examples relevant to us
  satisfy the condition that $H/V$ is
  solvable-by-finite.
\end{remark}
}

These two properties are important because of the following
result:
\begin{theorem}\label{theo:quotients_exist}\cite[Theorem
  3.46]{Linnell-Schick(2000)}:
  Let $H$ be a group which is cohomologically complete and
  which has enough amenable torsion-free quotients. Assume that there is a finite
  model for the classifying space of $H$.

  Let $1\to H\to G\to Q\to 1$ be an extension with finite quotient $Q$
  and such that $G$ is torsion-free. Then there is a normal subgroup
  $U$ of $G$, contained in $H$, such that $G/U$ is torsion-free and
  elementary amenable.
\end{theorem}
The proof of this theorem uses Sylow's theorems to reduce to the case
where $Q$ is a $p$-group. One then uses the pro-$p$ completions as an
intermediate
step to show that if no such torsion-free quotient exists, then the
projection $G\to Q$ induces a split injective map in
cohomology (we might have to replace $Q$ with a non-trivial
subgroup first). Using an Atiyah-Hirzebruch spectral sequence
argument, we
show that the same is true for stable cohomotopy. But then a
fixed-point theorem of Jackowski implies that
the map $G\to Q$ itself splits, which contradicts the assumption that
$G$ is torsion-free.

We now obtain our main technical result:
\begin{theorem}\label{theo:abstract_crit}
  Let $H$ be a group with finite classifying space which is
  cohomologically complete and which has enough amenable torsion-free
  quotients. Assume $H$ satisfies the Baum-Connes conjecture. Let $1\to H\to G\to Q\to 1$ be an extension of groups
  such that $G$ is torsion-free and such that $Q$ satisfies the
  Baum-Connes conjecture. Then the same is true for $G$.
\end{theorem}
\begin{proof}
  Because of Theorem \ref{theo:extensions} it suffices to prove the
  result if $Q$ is finite\longversion{ (and still $G$ torsion-free)}. 
 In this case,
  by Theorem \ref{theo:quotients_exist} we can replace the extension
  $1\to H\to G\to Q\to 1$ by another extension $1\to U\to G\to A\to 1$
  where $U$ is a subgroup of $H$ and $A$ is torsion-free and
  elementary amenable. Theorem \ref{theo:subgroups} implies that $U$
  fulfills the Baum-Connes conjecture. Because of
  Corollary \ref{corol:tfamext} the same is then true for $G$.
\end{proof}

In the last part of this note, we discuss examples of groups $H$, to
which Theorem \ref{theo:abstract_crit} applies.
\begin{definition}\label{def:groupclasses}
   Let $\extendableGroups$ be the class of groups $G$ which fulfill the
  following properties:
 $G$ has a finite classifying space,
 $G$ is cohomologically complete, and
 $G$ has enough nilpotent torsion-free quotients.

  Let $\BCclass$ be the class of groups defined in \cite[Definition
  5.22]{MR2027168}. 

  Set $\GoodextendableGroups:=\extendableGroups\cap\BCclass$.
\end{definition}

\begin{remark}\label{rem:prop_of_BCclass}
  We do not recall its complicated definition, but remark
  that it has the following important properties:
The class $\BCclass$ is closed under passing to subgroups, under extensions
with torsion free quotients and under finite products and finite free
  products. More generally, it is closed under passing to fundamental groups
  of graphs of groups. $\BCclass$ contains in particular all
one-relator groups and all Haken 3-manifold groups (and hence all knot or link
  groups), as well as all a-T-menable groups (in particular all amenable
  groups).
All these facts of the class $\BCclass$ and more information can be found in
  \cite{MR2027168}, in particular \cite[Theorem 5.23]{MR2027168}.  
\end{remark}

The following result is stated as \cite[Theorem 5.23]{MR2027168} 
or \cite[Theorem 5.2]{MR2181833}.
\begin{proposition}\label{prop:Utf_fulfills_BC}
  If $H\in\BCclass$ then the Baum-Connes
  conjecture with coefficients is true for $H$.
\end{proposition}

\begin{corollary}\label{corol:goodgroups}
  Assume that $H\in\GoodextendableGroups$ and we have an extension
  \begin{equation*}
1\to
  H\to G\to Q\to 1
\end{equation*}
  where $Q$ fulfills the Baum-Connes conjecture with coefficients and $G$ is
  torsion-free. Then the Baum-Connes conjecture with coefficients
   is true for $G$.
\end{corollary}
\begin{proof}
  Because $H\in\GoodextendableGroups\subset \BCclass$, by Proposition
  \ref{prop:Utf_fulfills_BC} $H$ satisfies the Baum-Connes conjecture. The
  definition of $\GoodextendableGroups$ implies that we can apply Theorem
  \ref{theo:abstract_crit} in our situation, which implies the assertion.
\end{proof}

\begin{definition}
    A semidirect product $H\semiprod Q$ of two groups $H$ and $Q$ is
  called \emph{$H_1$-trivial} if the induced action of $Q$ on
  $H_1(H,\integers)$ is trivial.
\end{definition}

\begin{definition}
  A one-relator group $G$ is called \emph{primitive} if it is finitely
  generated and if it has a presentation $G=\generate{x_1,\dots x_d \mid
  r}$ such that the element $r$ in the free group $F$ generated by
  $x_1,\dots,x_d$ is contained in the lower central series subgroup
  $\gamma_n(F)$ but not in
  $\gamma_{n+1}(F)$ and the image of $r$ in
  $\gamma_n(F)/\gamma_{n+1}(F)$ is not a proper power.
\end{definition}

\begin{example}\label{ex:primitive_one_relator}
Fundamental groups of orientable two-dimensional surfaces are
primitive one-relator groups, as well as
 one-relator groups where the least common multiple of the
    exponent-sums for the different generators $x_1,\dots,x_n$ in the relator
    $r$ is one.
\end{example}
\begin{proof} For the convenience of the reader, we sketch proofs of these
  well known facts.

  In the second case, the image of $r$ in the abelianization
  $\gamma_1(F)/\gamma_2(F)= F/[F,F]$ of the free group is a product of
  multiples $a_1[x_1], \dots, a_n[x_n]$ of the
  generators $[x_1],\dots,[x_n]$ of the free abelian group $F/[F,F]$ with
  generators $[x_1],\dots,[x_n]$ with least common multiple of
  $(a_1,\dots,a_n)$ equal to $1$, so is non-zero and not a proper power.

  The standard presentation $\generate{x_1,y_1,\dots,x_n,y_n\mid
    [x_1,y_1]\dots[x_n,y_n]}$ of a the fundamental group of an orientable
  surface of genus $n$ shows that the relation $r$, the product of the
  commutators $[x_i,y_i]$ is contained in $\gamma_2(F)$, but is a product of
  certain free generators of the free abelian group $\gamma_2(F)/\gamma_3(F)$,
  (this abelian group is freely generated by the images of all non-trivial
  commutators of the generators $x_1,\dots,x_n$, as follows from Magnus
  characterization of $\gamma_n(F)$, compare \cite{0016.29401}). Consequently,
  the image of $r$ in $\gamma_2(F)/\gamma_3(F)$ is non-zero and not a proper
  power. 
\end{proof}

\begin{definition} We denote the fundamental group $G$ of the
  complement of a tame link with $d$
  components in
  $S^3$  a \emph{link group with $d$ components}.
  We define the  \emph{linking diagram} to be the edge-labeled graph whose
  vertices
  are the components of the link, and such that any pair of vertices
  is joined by exactly one edge. Each edge is labeled with the linking
  number of the two link
  components involved. 

  We say the link group $G$ is \emph{primitive} if for each prime $p$
  there is a spanning subtree of the linking diagram such that none of
  the labels of the edges of this subtree is congruent to $0$ modulo
  $p$.

  Observe that in particular every \emph{knot group} (i.e.~a link group 
  with only one component) is primitive, as the linking diagram has one vertex
  and no edges, i.e.~a tree where we don't have to worry about any label of
  any dege.
\end{definition}

\begin{theorem}\label{theo:examples}
  The class $\GoodextendableGroups$ is closed under
  $H_1$-trivial semidirect products. It contains
  \begin{itemize}
  \item all primitive link groups
  \item  all primitive one-relator groups
  \item the fundamental groups of all fiber-type arrangements (as
    defined in \cite{Falk-Randell(1985)})
  \item Artin's pure braid groups $P_n$.
\end{itemize}
\end{theorem}
\begin{proof}
  It follows from \cite[Theorems 5.26, 5.40, Corollary
  5.27, Propositions 5.30, 5.34]{Linnell-Schick(2000)}, using
  \cite{Linnell-Schick(2000c)} and \cite{kuempel07:_towar_atiyah} that the
  groups which are mentioned belong to $\mathcal{F}$. By Remark
  \ref{rem:prop_of_BCclass} every link-group and every $1$-relator group
  belongs to $\BCclass$. Moreover, Artin's pure braid groups are iterated
  extensions of free groups. Since free groups belong to $\BCclass$ and
  $BCclass$ is closed under extensions with torsion-free quotients, the pure
  braid groups belong to $\BCclass$. Exactly the same argument applis to
  general fiber-type arrangements as defined in \cite{Falk-Randell(1985)}.
\end{proof}

\begin{corollary}
  All the torsion-free finite extensions of the groups listed in
  Theorem \ref{theo:examples}, in particular Artin's full braid
  groups $B_n$, satisfy the Baum-Connes conjecture.
\end{corollary}
\begin{proof}
  It is a classical result that the full braid group is torsion-free, and it
  is a finite extension of the pure braid group. Theorem \ref{theo:examples}
  implies that we can apply Corollary \ref{corol:goodgroups} to the extension
  \begin{equation*}
    1\to P_n\to B_n\to B_n/P_n\to 1,
  \end{equation*}
  or the corresponding general exact sequence, and the assertion follows.
\end{proof}

In Corollary \ref{corol:goodgroups} we can relax the condition that
$H$ belongs to $\GoodextendableGroups$ a little bit:

\begin{proposition}
    Assume $G_1,\dots,G_N\in \GoodextendableGroups$, and
    $Q\in\BCclass$ has a finite
  classifying space,
  is cohomologically complete, and has enough amenable torsion-free
  quotient. For example, $Q$ could be a free product $Q=Q_1*\cdots *
  Q_m$ of finitely many groups 
  $Q_i\in\GoodextendableGroups$. Define
  \begin{equation*}
    G:= G_1\semiprod\bigl(G_2\semiprod(\cdots (G_N\semiprod Q))\bigr),
  \end{equation*}
  where each semidirect product is $H_1$-trivial. Then $G$ has a
  finite classifying space, enough
 amenable torsion-free quotients, and is cohomologically complete. If there is an exact
  sequence 
  \begin{equation*}
    1\to G\to H\to A\to 1
  \end{equation*}
  with $H$ torsion-free such that $A$ satisfies the Baum-Connes
  conjecture, then the same is true for $H$.
\end{proposition}
\begin{proof}
  In \cite[Proposition 4.30]{Linnell-Schick(2000)} we prove that $G$
  has the desired properties. Because $\BCclass$ is closed under extension,
  $G\in \BCclass$, in particular $G$ fulfills the
  Baum-Connes conjecture. The last statement now
  follows from Theorem \ref{theo:abstract_crit}.
\end{proof}

{\small
\bibliographystyle{plain}
\bibliography{Baum_Connes_ext}
}
\end{document}